\documentclass{amsart}
\usepackage{amssymb, amsmath}
\usepackage[latin1]{inputenc}
\usepackage[final]{graphicx}
\newtheorem{theorem}{Theorem}[section]
\newtheorem{lemma}[theorem]{Lemma}
\newtheorem{remark}[theorem]{Remark}

\makeatletter
 \@addtoreset{equation}{section}
\makeatother

\begin{document}

\title{The affine quasi-Einstein Equation for homogeneous surfaces}
\author[M. Brozos-V\'{a}zquez et al.]
{M. Brozos-V\'{a}zquez \, E. Garc\'{i}a-R\'{i}o\, P. Gilkey,  and \, X. Valle-Regueiro}
\address{MBV: Universidade da Coru\~na, Differential Geometry and its Applications Research Group, Escola Polit\'ecnica Superior, 15403 Ferrol,  Spain}
\email{miguel.brozos.vazquez@udc.gal}
\address{EGR-XVR: Faculty of Mathematics,
University of Santiago de Compostela,
15782 Santiago de Compostela, Spain}
\email{eduardo.garcia.rio@usc.es, javier.valle@usc.es}
\address{PBG: Mathematics Department, University of Oregon, Eugene OR 97403-1222, USA}
\email{gilkey@uoregon.edu}
\thanks{Suported by Project MTM2016-75897-P (AEI/FEDER, UE)}
\begin{abstract}
We study the affine quasi-Einstein Equation for homogeneous surfaces. This gives rise through the modified Riemannian extension to new half conformally flat
generalized quasi-Einstein neutral signature $(2,2)$ manifolds, to conformally Einstein manifolds and also to new Einstein manifolds
through a warped product construction.
\end{abstract}

\maketitle

\section{Introduction}

The affine quasi-Einstein Equation (see Equation~(\ref{E1.c})) is a $0^{th}$ order perturbation of the Hessian. 
It is a natural linear differential equation in affine differential geometry. 
We showed (see \cite{BGGV17a}) that it gives rise to strong projective invariants of the affine structure.
Moreover, this equation also appears in the classification of half conformally flat quasi-Einstein manifolds in signature $(2,2)$. 
In this paper, we will examine the solution space to the affine quasi-Einstein Equation in the context of homogeneous affine geometries.

A description of locally homogeneous affine surfaces has been given by Opozda \cite{Op04} (see Theorem~\ref{T1.7} below). They fall into 3 families.
The first family is given by the Levi-Civita connection of a surface of constant curvature (Type $\mathcal{C}$). There are two
other families. The first (Type~$\mathcal{A}$) generalizes the Euclidean connection and the second (Type~$\mathcal{B}$) is a
generalization of the hyperbolic plane. As the Type $\mathcal{C}$ geometries are very rigid, we shall focus on the other two geometries. There
are many non-trivial solutions of the affine quasi-Einstein Equation for Type $\mathcal{A}$ geometries (see Section \ref{ss1.5}) and for Type $\mathcal{B}$ 
geometries (see Section \ref{ss1.6}).
This leads (see Theorem~\ref{T1.1} and Remark~\ref{R1.2}) to new examples of half conformally flat and conformally Einstein isotropic quasi-Einstein manifolds of signature $(2,2)$. We also use results of \cite{KimKim} to construct new higher dimensional Einstein manifolds.
Our present discussion illustrates many of the results of \cite{BGGV17a} and focusses on the dimension of
the eigenspaces of the solutions to the affine quasi-Einstein Equation for homogeneous surfaces.

\subsection{Notational conventions}
Recall that a pair $\mathcal{M}=(M,\nabla)$ is said to be an {\it affine manifold}
if $\nabla$ is a torsion free connection on the tangent bundle of a smooth manifold $M$ of dimension $m\ge2$. 
We shall be primarily interested in the case of affine surfaces ($m=2$) but it is convenient to work in greater generality
for the moment. In a system of local coordinates, express 
$\nabla_{\partial_{x^i}}\partial_{x^j}=\Gamma_{ij}{}^k\partial_{x^k}$ where we adopt the Einstein convention and sum over repeated indices. 
The connection $\nabla$ is torsion free if and only if the Christoffel symbols $\Gamma=(\Gamma_{ij}{}^k)$
satisfy the symmetry $\Gamma_{ij}{}^k=\Gamma_{ji}{}^k$ or, equivalently, if given any point $P$ of $M$, there
exists a coordinate system centered at $P$ so that in that coordinate system we have $\Gamma_{ij}{}^k(P)=0$. 

Let $f$ be a smooth function on $M$. The Hessian
\begin{equation}\label{E1.a}
\mathcal{H}_\nabla f=\nabla^2f:=(\partial_{x^i}\partial_{x^j}f-\Gamma_{ij}{}^k\partial_{x^k}f)\,dx^i\otimes dx^j
\end{equation}
is an invariantly defined symmetric $(0,2)$-tensor field; $\mathcal{H}_\nabla:C^\infty(M)\rightarrow C^\infty(S^2(M))$ is a second order partial differential operator
which is natural in the context of affine geometry.
The curvature operator $R_\nabla$ and the Ricci tensor $\rho_\nabla$ are defined by setting:
$$
R_\nabla(x,y):=\nabla_x\nabla_y-\nabla_y\nabla_x-\nabla_{[x,y]}\text{ and }
\rho_\nabla(x,y):=\operatorname{Tr}\{z\rightarrow R_\nabla(z,x)y\}\,.
$$ 
The Ricci tensor carries the geometry if $m=2$; an affine surface is flat if and only if $\rho_\nabla=0$ because
$$
\rho_{11}=R_{211}{}^2,\quad\rho_{12}=R_{212}{}^2,\quad\rho_{21}=R_{121}{}^1,\quad
\rho_{22}=R_{122}{}^1\,.
$$
In contrast to the situation in Riemannian geometry, $\rho_\nabla$ is not in general a symmetric $(0,2)$-tensor field.
The symmetrization and anti-symmetrization of the Ricci tensor are defined by setting, respectively,
$$\textstyle\rho_{s,\nabla}(x,y):=\frac12\{\rho_\nabla(x,y)+\rho_\nabla(y,x)\}\text{ and }
\textstyle\rho_{a,\nabla}(x,y):=\frac12\{\rho_\nabla(x,y)-\rho_\nabla(y,x)\}\,.
$$
We use $\rho_{s,\nabla}$ to define a $0^{\operatorname{th}}$
order perturbation of the Hessian. The {\it affine quasi-Einstein operator}
$ \mathfrak{Q}_{\mu,\nabla}:C^\infty(M)\rightarrow C^\infty(S^2(M))$ is defined by setting:
\begin{equation}\label{E1.b}
 \mathfrak{Q}_{\mu,\nabla} f:=\mathcal{H}_\nabla f-\mu f\rho_{s,\nabla}
 \,.
\end{equation}
The eigenvalue $\mu$ is a parameter of the theory; again, this operator is natural in the category of affine manifolds.
The {\it affine quasi-Einstein Equation} is the equation:
\begin{equation}\label{E1.c}
\mathfrak{Q}_{\mu,\nabla} f=0\text{ i.e. }\mathcal{H}_\nabla f=\mu f\rho_{s,\nabla}\,.
\end{equation}
We introduce the associated eigenspaces by setting:
$$
E(\mu,\nabla):=\ker( \mathfrak{Q}_{\mu,\nabla})=\{f\in C^\infty(M):\mathcal{H}_\nabla f=\mu f\rho_{s,\nabla}\}\,.
$$
Similarly, if $P$ is a point of $M$, we let $E(P,\mu,\nabla)$ be the space of germs of solutions to 
Equation~(\ref{E1.c}) which are defined near $P$. Note that $E(0,\nabla)=\ker(\mathcal{H}_\nabla)$
is the set of {\it Yamabe solitons}. Also note that $\rho_{s,\nabla}=0$ implies $E(\mu,\nabla)=E(0,\nabla)$ for any $\mu$.
If $\mu\ne0$ and $f>0$, let $\hat f:=-2\mu^{-1}\log(f)$, i.e. $f=e^{-\frac12\mu\hat f}$.
This transformation converts Equation~(\ref{E1.c}) into the equivalent non-linear equation:
\begin{equation}\label{E1.d}
\mathcal{H}_\nabla\hat f+2\rho_{s,\nabla}-\textstyle\frac12\mu d\hat f\otimes d\hat f=0\,.
\end{equation}

\subsection{Half conformally flat 4-dimensional geometry} Equation~(\ref{E1.d}) plays
an important role in the study of the quasi-Einstein Equation in neutral signature geometry \cite{BGGV17}.
Let $\mathcal{N}=(N,g,F,\mu_N)$ be a quadruple where $(N,g)$ is a 
pseudo-Riemannian manifold of dimension $n$, $F\in {C}^\infty(N)$,
and $\mu_N\in \mathbb{R}$. Let $\nabla^g$ be the Levi-Civita connection of $g$; the associated Ricci tensor $\rho_g$ is a symmetric $(0,2)$-tensor field.
We say that $\mathcal{N}$ is a \emph{quasi-Einstein manifold} if
$$
\mathcal{H}_{\nabla^g}F+\rho_g-\mu_NdF\otimes dF=\lambda\, g\text{ for }\lambda\in\mathbb{R}\,.
$$
We say $\mathcal{N}$ is {\it isotropic} if $\|dF\|=0$. We restrict to the 4-dimensional setting where Walker geometry 
(see \cite{DeR,Walker}) enters by means of the deformed Riemannian extension.
If $(x^1,x^2)$ are local coordinates on an affine surface $\mathcal{M}=(M,\nabla)$, 
let $(y_1,y_2)$ be the corresponding
dual coordinates on the cotangent bundle $T^*M$; if $\omega$ is a 1-form, then we can express $\omega=y_1dx^1+y_2dx^2$. Let $\Phi$ be an 
auxiliary symmetric $(0,2)$-tensor field.
The {\it deformed Riemannian extension} is defined \cite{CGGV09} by setting:
\begin{equation}\label{E1.e}
g_{\nabla,\Phi}= 2 dx^i \circ dy_i +\left\{-2y_k \Gamma_{ij}{}^k+\Phi_{ij}\right\}dx^i \circ dx^j\,.
\end{equation}
These neutral signature metrics are invariantly defined. Let $\pi:T^*M\rightarrow M$ be the natural projection. 
One has the following useful intertwining relation:
$$
\mathcal{H}_{g_{\nabla,\Phi}}\pi^*\hat f=\pi^*\mathcal{H}_\nabla\hat f,
\quad\rho_{g_{\nabla,\Phi}}=2\pi^*\rho_{s,\nabla},\quad\|d\pi^*\hat f\|_{g_{\nabla,\Phi}}
^2=0\,,
$$
for any $\hat f\in C^\infty(M)$.
The following observation is now immediate; note the factor of $\frac12$ in passing from the eigenvalue $\mu$ on the base to the eigenvalue 
$\mu_{T^*M}$ on the total  space:

\begin{theorem}\label{T1.1}
Let $(M,\nabla)$ be an affine surface 
and let $\hat f\in C^\infty(M)$ satisfy Equation~(\ref{E1.d}) or, equivalently, $f=e^{-\frac12\mu\hat f}\in E(\mu,\nabla)$.
Let $F=\pi^*\hat f$,
and let $\Phi$ be an arbitrary symmetric $(0,2)$-tensor field on $M$. Then
$(T^*M,g_{\nabla,\Phi},F,\mu_{T^*M})$ for $\mu_{T^*M}=\frac12\mu$ is a self-dual isotropic quasi-Einstein Walker manifold of signature $(2,2)$ with $\lambda=0$.
\end{theorem}

\begin{remark}\label{R1.2}\rm Starting with a quasi-Einstein manifold $(N, g, F, \mu_N)$ where $\mu_N=\frac{1}{r}$ for $r$ a positive integer, 
there exist appropriate Einstein fibers $E$ of dimension $r$ so that the warped product
$N\times_\varphi E$ is Einstein where $\varphi=e^{-F/r}$ \cite{KimKim}.
We can use Theorem~\ref{T1.1} to construct  self-dual isotropic quasi-Einstein Walker manifolds of neutral signature $(2,2)$ from
affine quasi-Einstein surfaces. Thus it is important to have solutions to the affine quasi-Einstein Equation for quite general $\mu$ and, in particular,
for $\mu=\frac 2r$. This will
be done quite explicitly in our subsequent analysis.
The parameter $\mu_N=-\frac{1}{n-2}$ is a distinguished value which is often exceptional. For $n\geq 3$, $(N,g, f,\mu_N= -\frac{1}{n-2})$ is a 
quasi-Einstein manifold if and only if $e^{-\frac{2}{n-2}f}g$ is Einstein \cite{EGR-KR2}. Therefore, taking into the account the fact that $\mu_{T^*M}=\frac12\mu_M$,
having solutions for the parameter $\mu_m=-\frac{1}{m-1}$ on $(M^m,\nabla)$ gives rise to conformally Einstein Riemannian extensions 
$(T^*M,g_{\nabla,\Phi})$ \cite{BGGV17}.
\end{remark}

The critical eigenvalue $\mu_m=-\frac1{m-1}$ is distinguished in this theory (see Theorem~\ref{T1.5} below). For affine surfaces, this corresponds
to $\mu_2=-1$ or, equivalently, $\mu_{T^*M}=-\frac12$. Excluding this value, we have the following converse to Theorem~\ref{T1.1} (see  \cite{BGGV17}). 

\begin{theorem}
Let $(N,g,F,\mu_N)$ be a self-dual quasi-Einstein manifold of signature $(2,2)$ which is not locally conformally flat with $\mu_N\neq -\frac{1}{2}$. Assume
$(N,g)$ is not Ricci flat. Then $\lambda=0$ and $(N,g,F,\mu_N)$ is locally isometric to a manifold which has the form given in
Theorem~\ref{T1.1}.
\end{theorem}

\subsection{Foundational results concerning the affine quasi-Einstein Equation}
 We established the following result in \cite{BGGV17a}:

\begin{theorem}\label{T1.4}
Let $\mathcal{M}=(M,\nabla)$ be an affine manifold. Let $f\in E(P,\mu,\nabla)$.
\begin{enumerate}
\item One has $f\in C^\infty(M)$. If $\mathcal{M}$ is real analytic, then $f$ is real analytic.
\item If $X$ is the germ of an affine Killing vector field based at $P$, then\newline$Xf\in E(P,\mu,\nabla)$.
\item If $f(P)=0$ and $df(P)=0$, then $f\equiv0$. Thus $\dim\{E(P,\mu,\nabla)\}\le m+1$.
\item If $M$ is simply connected and if $\dim\{E(P,\mu,\nabla)\}$ is constant on $M$, then $f$ extends
uniquely to an element of $E(\mu,\nabla)$.
\end{enumerate}
\end{theorem}

We say that  $\nabla$ and $\tilde\nabla$  are {\it projectively equivalent} if there exists a $1$-form
$\omega$ so that
$\nabla_XY= \tilde\nabla_XY+\omega(X)Y+\omega(Y)X$ for all $X$ and $Y$.
The equivalence is said to be a
{\it strong projective equivalence} if $\omega$ is closed.  If two projectively equivalent connections have symmetric Ricci tensors, then 
the two connections are, in fact, strongly projectively equivalent \cite{E70, NS, S95}.
A connection $\nabla$ is said to be {\it projectively flat} (resp. {\it strongly projectively flat}) if $\nabla$
is projectively equivalent (resp. strongly projectively equivalent) to a flat connection.

\begin{theorem}\label{T1.5}
If $\mathcal{M}$ is an affine surface, then $\dim\{E(\mu_2,\nabla)\}\ne2$. Moreover 
\begin{enumerate}
\item $\mathcal{M}$ is strongly projectively flat if and only if $\dim\{E(\mu_2,\nabla)\}=3$.
\item If $\mathcal{M}$ is strongly projectively flat and $\operatorname{rank}\rho_\nabla=2$, then for $\mu\neq \mu_2$:
\begin{enumerate}
	\item $\dim \{E(\mu,\nabla)\}=0$ for $\mu\neq 0$, and
	\item $\dim \{E(0,\nabla)\}=1$.
\end{enumerate}
\item If $\dim\{E(\mu,\nabla)\}=3$ for $\mu\ne\mu_2$, then $\mathcal{M}$ is Ricci flat and also strongly projectively flat.
\end{enumerate}
\end{theorem}

\begin{remark}\rm In Theorem \ref{T1.17}, we will exhibit
Type $\mathcal{B}$ surfaces with $\operatorname{rank}\rho_{s,\nabla}=2$ and $\operatorname{dim}\{ E(\mu,\nabla)\}\neq 0$ for
$\mu\ne0$. Consequently, Theorem~\ref{T1.5}~(2) fails without the assumption that $\mathcal{M}$ is strongly projectively flat.
\end{remark}

\subsection{Locally homogeneous surfaces}
Let $\Gamma_{ij}{}^k=\Gamma_{ji}{}^k\in\mathbb{R}$ define a connection $\nabla^\Gamma$ on $\mathbb{R}^2$.
The translation subgroup $(x^1,x^2)\rightarrow(x^1+a^1,x^2+a^2)$ acts transitively on $\mathbb{R}^2$ and preserves $\nabla^\Gamma$
so $(\mathbb{R}^2,\Gamma)$ is a homogeneous geometry.
In a similar fashion, let $\Gamma_{ij}{}^k=(x^1)^{-1}C_{ij}{}^k$ for $C_{ij}{}^k=C_{ji}{}^k\in\mathbb{R}$ define
a connection $\nabla^C$ on $\mathbb{R}^+\times\mathbb{R}$. The non-Abelian group $(x^1,x^2)\rightarrow(ax^1,ax^2+b)$ for $a>0$ and
$b\in\mathbb{R}$
acts transitively on the geometry $(\mathbb{R}^+\times\mathbb{R},\nabla^C)$
so this also is a homogeneous geometry. The following result was established by  Opozda \cite{Op04}.

\begin{theorem}\label{T1.7}
Let $\mathcal{M}=(M,\nabla)$ be a locally homogeneous affine surface. Then at least one of the following
three possibilities hold describing the local geometry:
\begin{itemize}
\item[($\mathcal{A}$)] There exists a coordinate atlas so the Christoffel symbols
$\Gamma_{ij}{}^k$ are constant.
\item[($\mathcal{B}$)] There exists a coordinate atlas so the Christoffel symbols have the form
$\Gamma_{ij}{}^k=(x^1)^{-1}C_{ij}{}^k$ for $C_{ij}{}^k$ constant and $x^1>0$.
\item[($\mathcal{C}$)] $\nabla$ is the Levi-Civita connection of a metric of constant Gauss
curvature.
\end{itemize}\end{theorem}

Let $\mathcal{M}=(M,\nabla)$ be a locally homogeneous affine surface which is not flat, i.e. $\rho_\nabla$ does not vanish identically.
One says that $\mathcal{M}$ is a {\it Type~$\mathcal{A}$}, {\it Type~$\mathcal{B}$} or
{\it Type~$\mathcal{C}$} surface depending on which possibility holds in Theorem~\ref{T1.7}. These are not exclusive
possibilities.
 
\begin{remark}
\rm
Type $\mathcal{C}$ surfaces are strongly projectively flat with Ricci tensor of rank $2$ in the non-flat case. Hence Theorem \ref{T1.5} shows that $\operatorname{dim}\{ E(0,\nabla)\}=1$, $\operatorname{dim}\{ E(-1,\nabla)\}=3$, and $\operatorname{dim}\{ E(\mu,\nabla)\}=0$ otherwise.
\end{remark}

We shall say that two Type~$\mathcal{A}$ surfaces are {\it linearly equivalent} if they are intertwined by the action of an element of the
general linear group. Similarly we shall say that two Type~$\mathcal{B}$ surfaces are {\it linearly equivalent} if they are intertwined
by some action $(x^1,x^2)\rightarrow(x^1,ax^1+bx^2)$ for $b\ne0$.
The Ricci tensor is symmetric for Type~$\mathcal{A}$ and Type~$\mathcal{C}$ geometries; there are Type~$\mathcal{B}$ geometries
where the Ricci tensor is not symmetric.

Let $\mu\in\mathbb{R}$.
It is convenient to consider complex solutions to the affine quasi-Einstein Equation
by taking $E_{\mathbb{C}}(\mu,\nabla):=E(\mu,\nabla)\otimes_{\mathbb{R}}\mathbb{C}$.
Real solutions can then be obtained by taking the real and imaginary parts as  both the underlying equation and the eigenvalue are real.

 \subsection{Type~$\mathcal{A}$ surfaces} \label{ss1.5}
The following result will be proved in Section~\ref{S3}; it is the foundation of our later results 
concerning Type~$\mathcal{A}$ geometry as it provides the ansatz for our
computations. Let $\mathfrak{A}$ be the commutative unital algebra of affine Killing vector fields generated by $\{\partial_{x^1},\partial_{x^2}\}$.

\begin{theorem}\label{T1.9}
Let $\mathcal{E}$ be a finite dimensional $\mathfrak{A}$ submodule of $C^\infty(\mathbb{R}^2)\otimes_{\mathbb{R}}\mathbb{C}$.
Then there exists a basis for $\mathcal{E}$ of functions of the form $e^{\alpha_1x^1+\alpha_2x^2}p(x^1,x^2)$ for
$p$ polynomial, where $\alpha_i\in\mathbb{C}$. Furthermore, $e^{\alpha_1x^1+\alpha_2x^2}\partial_{x^i}p\in\mathcal{E}$ for $i=1,2$.
\end{theorem}

Let $\mathcal{M}=(\mathbb{R}^2,\nabla)$ be a Type~$\mathcal{A}$ surface.
Any Type~$\mathcal{A}$ surface is strongly projectively flat with symmetric Ricci tensor \cite{BGG16}.
The following result will be established in Section~\ref{S4}.

\begin{theorem}\label{T1.10}
Let $\mathcal{M}$ be a Type~$\mathcal{A}$ surface which is not flat.
\begin{enumerate}
\item Let $\mu=0$. Then $E(0,\nabla)=\operatorname{Span}\{1\}$ or, up to linear equivalence, one of the following holds:
\begin{enumerate}
\item $\Gamma_{11}{}^1=1$, $\Gamma_{12}{}^1=0$, $\Gamma_{22}{}^1=0$, and $E(0,\nabla)=\operatorname{Span}\{1,e^{x^1}\}$.
\smallbreak\item $\Gamma_{11}{}^1=\Gamma_{12}{}^1=\Gamma_{22}{}^1=0$, and $E(0,\nabla)=\operatorname{Span}\{1,x^1\}$.
\end{enumerate}
\item Let $\mu=-1$. Then $\dim\{E(-1,\nabla)\}=3$.
\item Let $\mu\ne 0,-1$. Then $\dim\{E(\mu,\nabla)\}=\left\{\begin{array}{ll}
2\text{ if } \operatorname{rank}\rho_{\nabla}=1\\
0\text{ if } \operatorname{rank}\rho_{\nabla}=2
\end{array}\right\}$.
\end{enumerate}\end{theorem}

\begin{remark}\label{Re1.11}
\rm
Let $\mathcal{M}$ be a Type $\mathcal{A}$ surface which is not flat. By Remark~\ref{R1.2},
since $\operatorname{dim}\{ E(-1,\nabla)\}=3$, the corresponding 
Riemannian extension $\mathcal{N}:=(T^*M,g_{\nabla,\Phi})$ is conformally Einstein for any $\Phi$. 
Although in this instance $\mathcal{N}$ is locally conformally flat if $\Phi=0$ \cite{Afifi},
$\mathcal{N}$ is not locally conformally flat for generic $\Phi\neq 0$.
\end{remark}

\subsection{Type~$\mathcal{B}$ surfaces}\label{ss1.6}
Let $\mathfrak{B}$ be the unital non-commutative algebra generated by the vector fields $\partial_{x^2}$ and
 $X:=x^1\partial_{x^1}+x^2\partial_{x^2}$. The representation theory of the algebra $\mathfrak{B}$ is crucial
 to our investigation:
  
 \begin{theorem}\label{T1.12}
 Let $\mathcal{E}$ be a finite dimensional $\mathfrak{B}$ submodule of $C^\infty(\mathbb{R}^+\times\mathbb{R})\otimes_{\mathbb{R}}\mathbb{C}$.  \begin{enumerate}
\item If $f\in\mathcal{E}$, then
$ f=\sum_{\alpha,i,j}c_{\alpha,i,j}(x^1)^\alpha(\log(x^1))^i(x^2)^j$ where in this finite sum $c_{\alpha,i,j}\in\mathbb{C}$,
$\alpha\in\mathbb{C}$, and $i$ and $j$ are non-negative integers.
\item If $\dim\{\mathcal{E}\}=1$, then  $\mathcal{E}=\operatorname{Span}_{\mathbb{C}}\{(x^1)^a\}$ for some $a$.
\item If $\dim\{\mathcal{E}\}=2$, then one of the following possibilities hold:
\begin{enumerate}
\item $\mathcal{E}=\operatorname{Span}_{\mathbb{C}}\{(x^1)^\alpha,(x^1)^\alpha(c_1x^1+x^2)\}$.
\item $\mathcal{E}=\operatorname{Span}_{\mathbb{C}}\{(x^1)^\alpha,(x^1)^\alpha\log(x^1)\}$.
\item $\mathcal{E}=\operatorname{Span}_{\mathbb{C}}\{(x^1)^\alpha,(x^1)^\beta\}$ for $\alpha\ne\beta$.
\end{enumerate}
\end{enumerate}
 \end{theorem}
 
 Let $\mathcal{M}$ be a Type~$\mathcal{B}$ surface. Then $\partial_{x^2}$ and $X$ are Killing vector fields for $M$ and hence,
 by Theorem~\ref{T1.4}, $E(\mu,\nabla)$ is a finite dimensional $\mathfrak{B}$ module. The functions of Assertion~(1) are all
 real analytic; this is in accordance with Theorem~\ref{T1.4}. We will assume $\mathcal{M}$ is not Ricci flat and thus
 $\dim\{E(\mu,\nabla)\}\le2$ for $\mu\ne-1$ so Assertions~(2,3) apply. An appropriate analogue of Assertion~(1) holds in
 arbitary dimensions.

 Surfaces of Type~$\mathcal{A}$ are strongly projectively flat. Thus any Type~$\mathcal{B}$ surface which is also Type~$\mathcal{A}$ is strongly projectively flat
 (see Theorem~\ref{T6.1}).
There are, however, strongly projectively flat surfaces of Type~$\mathcal{B}$ which are not of Type~$\mathcal{A}$. We will establish the following result
in Section~\ref{S6.1}.

\begin{theorem}\label{T1.13}
If $\mathcal{M}$ be a Type~$\mathcal{B}$ surface, then $\mathcal{M}$ is strongly projectively flat if and only if 
$\mathcal{M}$ is linearly equivalent to one of the surfaces:
\begin{enumerate}
\item $C_{12}{}^1=C_{22}{}^1=C_{22}{}^2=0$ (i.e. $\mathcal{M}$ is also of Type~$\mathcal{A}$).
\item $C_{11}{}^1=1+2v$, $C_{11}{}^2=0$, $C  _{12}{}^1=0$, $C_{12}{}^2=v$, $C_{22}{}^1=\pm1$, $C_{22}{}^2=0$.
\end{enumerate}
\end{theorem}

We remark that the special choice of $v=-1$ in Assertion~(2) corresponds to the hyperbolic plane and the Lorentzian analogue in Theorem~\ref{T6.1}~(3). 

We now turn to the study of the affine quasi-Einstein Equation.
We first examine the Yamabe solitons, working modulo linear equivalence. We shall establish the following result in Section~\ref{S6.2}.

\begin{theorem}\label{T1.14}
Let $\mathcal{M}$ be a Type~$\mathcal{B}$ surface. Then $E(0,\nabla)=\operatorname{Span}\{1\}$ except in the following cases where we
also require $\rho_\nabla\ne0$.
\begin{enumerate}
\item $C_{11}{}^2=c\,C_{11}{}^1$, $C_{12}{}^2=c\,C_{12}{}^1$, $C_{22}{}^2=c\,C_{22}{}^{ 1}$, $E(0,\nabla)=\operatorname{Span}\{1,x^2-cx^1\}$.
\item $C_{11}{}^1=-1$, $C_{12}{}^1=0$, $C_{22}{}^1=0$, $E(0,\nabla)=\operatorname{Span}\{1,\log(x^1)\}$.
\item $C_{11}{}^1=\alpha-1$, $C_{12}{}^1=0$, $C_{22}{}^1=0$, $E(0,\nabla)=\operatorname{Span}\{1,(x^1)^\alpha\}$ for $\alpha\ne0$.
\end{enumerate}
\end{theorem}

By Theorem~\ref{T1.5}~(1), $\dim\{E(-1,\nabla)\}=3$ if and only if $\mathcal{M}$ is strongly projectively flat. 
Thus this case is covered by Theorem~\ref{T1.13}.  Furthermore, $\dim\{E(-1,\nabla)\}\ne2$ by Theorem~\ref{T1.5}. 
We now examine the remaining case where $\dim\{E(-1,\nabla)\}=1$.

\begin{theorem}\label{T1.15}
Let $\mathcal{M}$ be a Type~$\mathcal{B}$ surface with $\rho_\nabla\ne0$ which is not strongly projectively flat. Then
$\dim\{E(-1,\nabla)\}=0$ except in the following cases where $\dim\{E(-1,\nabla)\}=1$:
\begin{enumerate}
\item $C_{22}{}^1=0$, $C_{22}{}^2=C_{12}{}^1\neq 0$.
\item $C_{22}{}^1=\pm 1$, $C_{12}{}^1=0$, $C_{22}{}^2=\pm 2C_{11}{}^2\neq 0$,
$C_{11}{}^1=1+2C_{12}{}^2\pm(C_{11}{}^2)^2$.
\end{enumerate} 
\end{theorem}

\begin{remark}\label{Re1.16}
\rm
Let $\mathcal{M}$ be a Type $\mathcal{B}$ surface with $\operatorname{dim}\{ E(-1,\nabla)\}=1$. 
By Remark~\ref{R1.2}, $\mathcal{N}:=(T^*M,g_{\nabla,\Phi})$ is conformally Einstein for any $\Phi$. Moreover, since 
$\mathcal{M}$ is not strongly projectively flat, unlike in the Type~$\mathcal{A}$
setting, $\mathcal{N}$ is not locally conformally flat for any $\Phi$ \cite{Afifi}.
\end{remark}

Let $\mu\ne0$ and $\mu\ne-1$. In the Type~$\mathcal{A}$ setting, Theorem \ref{T1.10} shows that $\dim\{E(\mu,\nabla)\}=0$ or $\dim\{E(\mu,\nabla)\}=2$.
The situation is quite different in the Type~$\mathcal{B}$ setting as there are examples where $\dim\{E(\mu,\nabla)\}=1$.

\begin{theorem}\label{T1.17} 
Let $\mathcal{M}$ be a Type~$\mathcal{B}$ surface which is not of Type~$\mathcal{A}$ with $\rho_{s,\nabla}\ne 0$. Let $\mu\notin\{0,-1\}$.
Then $\dim\{E(\mu,\nabla)\}\geq 1$ if and only if $\mathcal{M}$ is linearly equivalent to a surface with
\medbreak\qquad\qquad
 $C_{12}{}^1=0$, $C_{22}{}^1=\pm 1$, $C_{22}{}^2=\pm 2C_{11}{}^2$,
\smallbreak\qquad\qquad$\mu=\frac{-(C_{11}{}^1)^2+2 C_{11}{}^1 C_{12}{}^2 \pm 2 (C_{11}{}^2)^2 -(C_{12}{}^2)^2+2 
C_{12}{}^2+1}{(C_{11}{}^1-C_{12}{}^2-1)^2}$,
$C_{11}{}^1-C_{12}{}^2-1\ne0$.
\medbreak\noindent Furthermore,
$\dim\{E(\mu,\nabla)\}=2$ if and only if $\mathcal{M}$ is
linearly equivalent to one of: 
\medbreak\qquad
{\rm(1)} $C_{11}{}^1=-1+C_{12}{}^2$\,,
$C_{11}{}^2=0$\,,
$C_{12}{}^1=0$\,,
$C_{22}{}^1=\pm 1$\,,
$C_{22}{}^2=0$,
\medbreak\qquad\qquad $\mu=\frac{1}{2}C_{12}{}^2\neq 0$.
\medbreak\qquad
{\rm(2)} $ C_{11}{}^1=-\frac{1}{2}(5\pm 16 (C_{11}{}^2)^2)$\,,
$C_{12}{}^1=0$\,,
$C_{12}{}^2=-\frac{1}{2}(3\pm 8(C_{11}{}^2)^2)$,
\medbreak\qquad\qquad
$C_{22}{}^1=\pm 1$\,,
$C_{22}{}^2=\pm 2 C_{11}{}^2$,
$\mu=-\frac{3\pm 8(C_{11}{}^2)^2}{4\pm 8(C_{11}{}^2)^2}$,
$C_{11}{}^2\notin\{ 0,\pm\frac{1}{\sqrt{2}}\}$\,.
\end{theorem}

\section{The proof of Theorem~\ref{T1.5}}

Assertions~(1) and (3) of Theorem~\ref{T1.5} follow from more general results of \cite{BGGV17a}.  The following result (see \cite{E70, NS})
characterizes strongly projectively flat surfaces and will be used in the proof of Theorem~\ref{T1.5}~(2).

\begin{lemma}\label{L2.1}
If $\mathcal{M}=(M,\nabla)$ is an affine surface, then the following assertions are equivalent:
\begin{enumerate}
\item $\mathcal{M}$ is strongly projectively flat.
\item $\rho_\nabla$ and $\nabla\rho_\nabla$ are totally symmetric.
\item $\rho_\nabla$ is symmetric and $\mathcal{M}$ is projectively flat.
\end{enumerate}
\end{lemma}

\subsection{The proof of Theorem~\ref{T1.5}~(2)} Let $\mathcal{M}$ be strongly projectively flat with $\operatorname{rank}\rho_\nabla=2$.
By Lemma~\ref{L2.1}, $\rho_{\nabla}=\rho_{s,\nabla}$.
We may covariantly differenciate the quasi-Einstein  Equation \eqref{E1.c} with respect to $\partial_{x^i}$ in local coordinates to get
\[
\mathcal{H}_{\nabla,jk;i} f=\mu \{(\partial_{x^i} f) \rho_{\nabla,jk}+ f \rho_{\nabla,jk;i}\}.
\]
We interchange the indices $i$ and $j$ and subtract to get
\[
\begin{array}{rcl}
R_{\nabla,ijk}{}^l (\partial_{x^l}f)&=& \mathcal{H}_{\nabla,ik;j} f-\mathcal{H}_{\nabla,jk;i} f\\
\noalign{\medskip}
&=&\mu \{(\partial_{x^j} f) \rho_{\nabla,ik}+ f \rho_{\nabla,ik;j}\}-\mu \{(\partial_{x^i} f) \rho_{\nabla,jk}+ f \rho_{\nabla,jk;i}\}.
\end{array}
\]
By Lemma~\ref{L2.1}, $\nabla\rho_\nabla$ is totally symmetric and the previous expression simplifies:
\[
R_{\nabla,ijk}{}^l (\partial_{x^l}f)=\mu \{(\partial_{x^j} f) \rho_{\nabla,ik}-(\partial_{x^i} f) \rho_{\nabla,jk}\}.
\]
Since $M$ is 2-dimensional the only curvature term is 
$$
R_{\nabla,12}:\partial_{x^i}\rightarrow\rho_{\nabla,2i} \partial_{x^1}-\rho_{\nabla,1i} \partial_{x^2} \,.
$$
Consequently, $(R_{\nabla,12} \partial_{x^i})f=-\mu (R_{\nabla,12} \partial_{x^i})f$ 
and, hence,
\[
(\mu+1) (R_{\nabla,12} \partial_{x^i})f=0.
\] 
This is a homogeneous system of linear equations. Because $\rho_\nabla$ has rank 2,  $R_{\nabla,12}$ has rank 2. If $\mu\neq -1,0$, then the only solutions are 
$\partial_{x^i}f=0$ and $f$ is necessarily constant. This shows that $\operatorname{dim}\{ E(\mu,\nabla)\}=0$ for all $\mu\neq 0,-1$. Moreover, if $\mu=0$ one has $R_{\nabla,12} \partial_{x^i}f=0$ and $f$ is constant, thus showing that $\operatorname{dim}\{ E(0,\nabla)\}=1$.
\qed

\section{The proof of Theorem~\ref{T1.9}}\label{S3}
Let $\mathfrak{A}$ be the commutative unital algebra generated by $\{\partial_{x^1},\partial_{x^2}\}$.
Let $\mathcal{E}$ be a finite dimensional $\mathfrak{A}$ submodule of $C^\infty(\mathbb{R}^2)\otimes_{\mathbb{R}}\mathbb{C}$.
If $f\in\mathcal{E}$, set
$$
\mathcal{E}(f):=\operatorname{Span}\{f,\partial_{x^1}f,\dots,\partial^k_{x^1}f,\dots\}\subset\mathcal{E}\,.
$$
As $\mathcal{E}$ is finite dimensional, there is a minimal dependence relation:
\begin{equation}\label{E3.a}
\textstyle\prod_i(\partial_{x^1}-\lambda_i)^{n_i}f=0\text{ for }\lambda_i\in\mathbb{C}\text{ distinct}\,.
\end{equation}
Let $f_i:=\prod_{j\ne i}(\partial_{x^1}-\lambda_j)^{n_j}f\in\mathcal{E}(f)$. If we fix $x^2$, then
$f_i(x^1,x^2)$ satisfies the constant coefficient ODE $(\partial_{x^1}-\lambda_i)^{n_i}f_i(x^1,x^2)=0$ with suitably
chosen initial conditions determined by $\{f_i(0,x^2),\partial_{x^1}f_i(0,x^2),\dots\}$. Consequently, we may express
$f_i(x^1,x^2)=\{(x^1)^{n_i-1} h _{n_i-1}(x^2)+\dots+ h _0(x^2)\}e^{\lambda_ix^1}$. Since Equation~(\ref{E3.a}) is to be a minimal dependence
relation for $f$, $ h _{n_i-1}$ does not vanish identically. Consequently, the functions $\{(\partial_{x^1}-\lambda_i)^kf_i\}_{0\le k\le n_i-1}$ for fixed $i$
are linearly independent. The different exponential terms do not interact and thus for dimensional reasons, the collection of all these functions 
forms a basis for $\mathcal{E}(f)$. Thus, $f$ can be
expressed in terms of these elements, i.e. any element of $\mathcal{E}$ can be expressed as a sum of functions
of the form $e^{\alpha x^1} \sum_{i} (x^1)^{i} h_{i}(x^2)$ where $\alpha\in \mathbb{C}$. A similar analysis of the $x^2$ dependence shows that we can express 
$f$ in the given form. 
We complete the proof by noting that we can express
\medbreak\hfill $e^{\alpha_1x^1+\alpha_2x^2}\partial_{x^1}p=(\partial_{x^1}-\alpha_1)\{e^{\alpha_1x^1+\alpha_2x^2}p(\vec x)\}\in\mathcal{E}$.
\hfill\vphantom{.}\qed

\section{Type~$\mathcal{A}$ surfaces}\label{S4}
We shall assume that $\rho\ne0$. By Theorem~\ref{T1.4}, $\dim\{E(\mu,\nabla)\}\le3$ and by Theorem~\ref{T1.5}, $\dim\{E(\mu,\nabla)\}\ne3$ if $\mu\ne\mu_2=-1$.
We use these facts implicitly in what follows.

\subsection*{The proof of Theorem~\ref{T1.10}~(1)}
Let $\mathcal{M}$ be a Type~$\mathcal{A}$ surface. It is immediate from the definition that $1\in E(0,\nabla)$.
Suppose there exists a non-constant function $f\in E(0,\nabla)$, i.e. $\dim\{E(0,\nabla)\}=2$. 
By Lemma~4.1 of \cite{BGG16}, $R_{\nabla,12}(df)=0$.
 This implies that $df$ belongs to the kernel of the curvature operator. Consequently after a suitable linear change of coordinates, we have
 $f=f(x^1)$ for any $f\in E(0,\nabla)$. We apply Theorem~\ref{T1.9} to see that we may assume $f(x^1,x^2)=p(x^1)e^{a_1x^1}$.
\smallbreak\subsection*{Case 1a.} Suppose that $a_1\ne0$. By Theorem~\ref{T1.9}, we may assume
$f(x^1)=e^{a_1x^1}$. Since $\mathcal{H}_\nabla f=0$,
 Equation~(\ref{E1.a}) implies $\Gamma_{12}{}^1=0$, $\Gamma_{22}{}^1=0$, and $\Gamma_{11}{}^1=a$. Thus $a$ is real.
Rescale the first coordinate to assume $\Gamma_{11}{}^1=1$. A direct computation shows $E(0,\nabla))=\operatorname{Span}\{1,e^{x^1}\}$
 which is the possibility of Assertion~(1a).
 \smallbreak\subsection*{Case 2a.} Suppose that $a_1=0$ so $f(x^1,x^2)=p(x^1)$ is a non-constant polynomial. We apply Theorem~\ref{T1.9}
 to assume $p$ is linear. Subtracting the constant term then permits us to assume $p(x^1)=x^1$. We then obtain
 $\Gamma_{11}{}^1=\Gamma_{12}{}^1=\Gamma_{22}{}^1=0$.  A direct computation shows $E(0,\nabla))=\operatorname{Span}\{1,x^1\}$ which
 is the possibility of Assertion~(1b). 
 \hfill\qed
 
 \subsection{The proof of Theorem~\ref{T1.10}~(2)} Results of \cite{BGG16} show that any Type~$\mathcal{A}$ geometry is
 strongly projectively flat. Theorem~\ref{T1.5} then shows $\dim\{E(-1,\nabla)\}=3$.
 
 \subsection{The proof of Theorem~\ref{T1.10}~(3)} Assume that $\mu\notin\{0,-1\}$.
 
\subsection*{Case 3a.} Let $\mathcal{M}$ be a Type~$\mathcal{A}$ surface with $\operatorname{rank}\rho_\nabla=1$.
As $\dim\{E(\mu,\nabla)\}\le2$,
it suffices to show $\dim\{E(\mu,\nabla)\}\ge2$.
We make a linear change of coordinates to assume $\rho_{{ \nabla},11}=\rho_{\nabla,12}=0$. By Lemma~2.3 of \cite{BGG16}, this
implies $\Gamma_{11}{}^2=\Gamma_{12}{}^2=0$. The affine quasi-Einstein Equation for $f(x^1,x^2)=e^{a_2x^2}$ reduces to the single equation:
$$
a_2^2-a_2\Gamma_{22}{}^2-\mu\rho_{22}=0\,.
$$
Generically, this has two distinct complex solutions which completes the proof in this special case. However,
we must deal with the case in which the discriminant $(\Gamma_{22}{}^2)^2+4\mu\rho_{22}=0$. Since $\mu\ne0$ and $\rho_{22}\ne0$, 
$\Gamma_{22}{}^2\ne0$. Thus there is a single exceptional value $\mu=-(\Gamma_{22}{}^2)^2/(4\rho_{22})$. The
affine quasi-Einstein Equation for $f(x^1,x^2)=x^2e^{a_2x^2}$ again reduces to a single equation
$$
(2a_2-\Gamma_{22}{}^2)(4+2a_2x^2-x^2\Gamma_{22}{}^2)=0\,.
$$
Let $a_2=\frac12\Gamma_{22}{}^2$ to ensure $x^2e^{a_2x^2}\in E(\mu,\nabla)$
so $\dim\{E(\mu,\nabla)\}\ge2$ as desired.
\subsection*{Case 3b} Suppose $\mathcal{M}$ is a Type~$\mathcal{A}$ surface with $\operatorname{rank}\rho_\nabla=2$. 
We apply Theorem~\ref{T1.5}~(2) to show that $\dim\{E(\mu,\nabla)\}=0$.\qed
\section{The proof of Theorem~\ref{T1.12}}
We prove the three assertions seriatim.

\subsection*{The proof of Theorem~\ref{T1.12}~(1)}
Let $\mathfrak{B}$ be the unital non-commutative algebra generated by $X:=x^1\partial_{x^1}+x^2\partial_{x^2}$ and $\partial_{x^2}$.
Let $\mathcal{E}$ be a finite dimensional $\mathfrak{B}$ submodule of $C^\infty(\mathbb{R}^+\times\mathbb{R})\otimes_{\mathbb{R}}\mathbb{C}$.
Let $0\ne f\in\mathcal{E}$. Applying exactly the same argument as that used to prove Theorem~\ref{T1.9}
permits us to expand $f$ in the form
$$
f=\sum_{ij}e^{\beta_jx^2}(x^2)^i h _{ij}(x^1)\text{ where }\beta_j\in\mathbb{C}\,.
$$
Suppose that $h_{ij}\not\equiv0$ for some $\beta_j\ne0$. Choose $i_0$ maximal so $h_{i_0j}\not\equiv0$.   We compute
$$
X^nf=e^{\beta_jx^2}\{\beta_j^n(x^2)^{i_0+n}h_{i_0j}(x^1)+O((x^2)^{i_0+n-1})\}+\dots
$$
where we have omitted terms not involving the exponential $e^{\beta_jx^2}$.
This constructs an infinite sequence of linearly independent functions in $\mathcal{E}$ which contradicts our assumption that
$\mathcal{E}$ is finite dimensional. Consequently $f$ is polynomial in $x^2$. We wish to determine the form of the
coefficient functions $h_i(x^1)$. Let $\tilde X:=x^1\partial_{x^1}$. We have
$X^kf=\sum_i(x^2)^i(i+\tilde X)^kh_i$. Since $\mathcal{E}$ is finite dimensional, the collection $\{(i+\tilde X)^kh_i\}$,
or equivalently $\{\tilde X^kh_i\}$, can not be an infinite sequence of linearly independent functions. If $h_i$ is constant,
we do not need to proceed further in considering $h_i$. Otherwise, there must be a minimal dependence relation
which we can factor in the form
\begin{equation}\label{E5.a}
\textstyle\prod_j(\tilde X-\lambda_j)^{n_j}h_i=0\,.
\end{equation}
The solutions to Equation~(\ref{E5.a})
are spanned by the functions $ (\log(x^1))^k(x^1)^\lambda$ but apart from that, the analysis is the same as
that performed in the proof of Theorem~\ref{T1.9} and we can expand each $h_i$ in terms of these functions. This completes
the proof of Assertion~(1).

\subsection*{The proof of Theorem~\ref{T1.12}~(2)} Let $\dim\{\mathcal{E}\}=1$ and let $0\ne f\in\mathcal{E}$. 
Expand $f$ in the form of Assertion~(1) and choose $j_0$ maximal so $c_{a,i,j_0}\ne0$ for some $(a,i)$. If $j_0>0$, then 
$\partial_{x^2}f\ne0$. Consequently $\{f,\partial_{x^2}f\}$
are linearly independent elements of $\mathcal{E}$ which is false. Let $c_{a,i}:=c_{a,i,0}$. We have
$$
\displaystyle Xf=\sum_{a,i}c_{a,i}(x^1)^a\{a\log(x^1)^i+i\log(x^1)^{i-1}\}\in\mathcal{E}\,.
$$
Since $\dim\{\mathcal{E}\}=1$, $Xf$ must be a multiple of $f$. 
Thus, there is exactly one value of $a$ so $c_{a,i_0}\ne0$. Furthermore, one has
$i_0=0$. This implies $f=(x^1)^a$ as desired.

\subsection*{The proof of Theorem~\ref{T1.12}~(3)}
 Let $\dim\{\mathcal{E}\}=2$. If Assertion~(3c) fails, we can choose
$0\ne f\in\mathcal{E}$ so that $f\ne(x^1)^a$ for any $a$. Expand $f$ in the form of Assertion~(1) and 
choose $j_0$ maximal so $c_{a,i,j_0}\ne0$. 
\smallbreak\noindent{\bf Step 3a.} Suppose $j_0\ge1$. Then $\{f,\partial_{x^2}f\}$ are linearly independent
elements of $\mathcal{E}$ and hence are a basis for $\mathcal{E}$. Consequently $\partial_{x^2}^2f=0$ and $f$
is linear in $x^2$. Let $\mathcal{E}_1:=\partial_{x^2}f\cdot\mathbb{C}$. This subspace is clearly invariant under
$X$ and $\partial_{x^2}$. 
Thus, by Assertion~(2),
$\partial_{x^2}f=(x^1)^\alpha$ for some $\alpha$. Consequently
$$
f(x^1,x^2)= h (x^1)+(x^1)^\alpha x^2\text{ for } h (x^1)=\sum_{a,i}c_{a,i}(x^1)^a(\log(x^1))^i\,.
$$
We have $(X-\alpha-1)f=(X-\alpha-1) h $ is independent of $x^2$ and thus belongs to $\mathcal{E}_1$. Consequently, this is a multiple of $(x^1)^\alpha$, i.e.
$$
\sum_{a,i}c_{a,i}(x^1)^a\{(a-\alpha-1)(\log(x^1))^i+i(\log(x^1))^{i-1}\}=c_0(x^1)^\alpha\,.
$$
We must therefore have $c_{a,i}=0$ for $a\notin\{\alpha+1,\alpha\}$.
If $a=\alpha+1$ or $a=\alpha$, we conclude $c_{\alpha,i}=0$ for $i>0$. We can eliminate the term involving $ (x^1)^\alpha $ by subtracting
an appropriate multiple of $\partial_{x^2}f= (x^1)^\alpha $. Thus $\mathcal{E}$ has the form given in Assertion~(3a).

\subsection*{Step 3b} Suppose $f=\sum_{a,i}c_{a,i}(x^1)^a(\log(x^1))^i$ is independent of $x^2$. We have
$$
(X-b)f=\sum_{a,i}c_{a,i}(x^1)^a\{(a-b)(\log(x^1))^i+i(\log(x^1))^{i-1}\}\,.
$$
Choose $i_0$ maximal so $c_{a,i_0}\ne0$. If $i_0=0$, then $f=\sum_a c_a(x^1)^a$. Since $f\ne(x^1)^a$, there are at least two
different exponents where $c_{a_i}\ne0$. Since $(X-a_i)f$ will be non-zero and not involve the exponent $a_i$, we conclude
for dimensional reasons that there are exactly two such indices and that $\mathcal{E}=\operatorname{Span}\{(x^1)^{a_1},(x^1)^{a_2}\}$
contrary to our assumption. Thus $i_0>0$ and $c_{a_1,i_0}\ne0$. Suppose there is more than 1 exponent. Then $\{f,(X-a_1)f,(X-a_2)f\}$
would be linearly independent. Thus we could take
$f=\sum_i(x^1)^a(\log(x^1))^i$. If $i_0\ge2$, we conclude $\{f,(X-a_1)f,(X-a_1)^2f\}$ are linearly
independent. Thus $f=(x^1)^a\{c_0+c_1 \log(x^1)\}$. Again, applying $(X-a)$, we conclude $(x^1)^a\in \mathcal{E}$ and thus
Assertion~(3b) holds.
\qed

\section{Type~$\mathcal{B}$ surfaces}

We refer to \cite{BGG16} for the proof of the following result.

\begin{theorem}\label{T6.1}
\ \begin{enumerate}
\item There are no surfaces which are both Type~$\mathcal{A}$ and Type~$\mathcal{C}$.
\item A type~$\mathcal{B}$ surface is locally isomorphic to a Type~$\mathcal{A}$ surface if and only if
$(C_{12}{}^1,C_{22}{}^1,C_{22}{}^2)=(0,0,0)$;
the Ricci tensor has rank $1$ in this instance.
\item A Type~$\mathcal{B}$ surface is locally isomorphic to a Type~$\mathcal{C}$ surface if and only
if $C$ is linearly equivalent to one of the following two examples whose non-zero Christoffel symbols are
$C_{11}{}^1=-1$, $C_{12}{}^2=-1$, $C_{22}{}^1=\pm1$.
The associated Type~$\mathcal{C}$ geometry is given by
$ds^2=(x^1)^{-2}\{(dx^1)^2\pm(dx^2)^2\}$.
\end{enumerate}
\end{theorem}

Throughout this section, we will let $\mathcal{M}$ be a Type~$\mathcal{B}$ affine surface; we assume $\rho\ne0$ to ensure the
geometry is not flat. In Section~\ref{S6.1}, we examine when a Type~$\mathcal{B}$ surface is strongly projectively flat. 
In Section~\ref{S6.2},
we determine the dimension of the space of Yamabe solitons $E(0,\nabla)$. 
In Section~\ref{S6.3}, we examine $\dim\{E(-1,\nabla)\}$.
In Section~\ref{S6.4}, we examine the general case where $\mu\notin\{0,-1\}$.

\subsection{The proof of Theorem~\ref{T1.13}}\label{S6.1}
By Lemma~\ref{L2.1}, $\mathcal{M}$ is strongly projectively flat if and only if $\rho_\nabla$ is symmetric and $\nabla\rho_\nabla$ is totally symmetric.
A direct calculation shows the surfaces in Assertions~(1) and (2) satisfy this condition. Conversely, $\rho_\nabla$ is symmetric if and only if $C_{22}{}^2=-C_{12}{}^1$.
Impose this condition. It is then immediate that $\nabla\rho_\nabla(1,2;1)=\nabla\rho_\nabla(2,1;1)$ and $\nabla\rho_\nabla(1,2;2)=\nabla\rho_\nabla(2,1;2)$.
The remaining two equations yield:
\begin{eqnarray*}
0&=&\nabla\rho_\nabla(1,2;1)-\nabla\rho_\nabla(1,1;2)\\
&=&C_{11}{}^1 C_{12}{}^1+3 C_{11}{}^2 C_{22}{}^1-2 C_{12}{}^1 C_{12}{}^2+2 C_{12}{}^1,\\
0&=&\nabla\rho_\nabla(1,2;2)-\nabla\rho_\nabla(2,2;1)\\
&=&2 C_{11}{}^1 C_{22}{}^1-6 (C_{12}{}^1)^2-4 C_{12}{}^2 C_{22}{}^1-2 C_{22}{}^1.
\end{eqnarray*}
Suppose first $C_{22}{}^1=0$. The second constraint yields $-6(C_{12}{}^1)^2=0$. Thus $C_{12}{}^1=0$ and $C_{22}{}^2=0$. This yields the
surfaces of Assertion~(1). We therefore assume $C_{22}{}^1\ne0$. We can rescale so $C_{22}{}^1=\varepsilon=\pm1$. 
Let $\tilde x^1=x^1$ and $\tilde x^2=c x^1+x^2$ define a shear. We then obtain
$\tilde C_{12}{}^1=C_{12}{}^1-cC_{22}{}^1$. Thus by choosing $c$ appropriately, we assume that $C_{12}{}^1=0$. We impose these constraints. Our equations
become, after dividing by { $\varepsilon$,} $3C_{11}{}^2=0$ and $-2+2C_{11}{}^1-4C_{12}{}^2=0$.
We set $C_{12}{}^2=v$ to obtain $C_{11}{}^1=1+2v$ and obtain the surfaces of Assertion~(2).\qed

\subsection{The proof of Theorem~\ref{T1.14}: $\mu=0$}\label{S6.2}
Let $\mathcal{M}$ be a Type~$\mathcal{B}$ surface which is not flat, or equivalently, not Ricci flat. Consequently,
$\dim\{E(0,\nabla)\}<3$ by Theorem~\ref{T1.5}. We have $1\in E(0,\nabla))$. Suppose, exceptionally, $\dim\{E(0,\nabla)\}=2$.
We examine the 3 cases of Theorem~\ref{T1.12} seriatim.
\subsection*{Case 1} Suppose $E(0,\nabla)=\operatorname{Span}\{(x^1)^\alpha,(x^1)^\alpha(-cx^1+x^2)\}$.
Since $1\in E(0,\nabla)$, we have $\alpha=0$ so $f=-cx^1+x^2$. The following equations yield Assertion~(1):
$$
(Q_{11}):\ 0=cC_{11}{}^1-C_{11}{}^2,\quad (Q_{12}):\ 0=cC_{12}{}^1-C_{12}{}^2,\quad (Q_{22}):\ 0=cC_{22}{}^1-C_{22}{}^2\,.
$$

\subsection*{Case 2} Suppose $E(0,\nabla)=\operatorname{Span}\{(x^1)^\alpha,(x^1)^\alpha\log(x^1)\}$.
Since $1\in E(0,\nabla)$, we have $\alpha=0$ so 
$f=\log(x^1)$. The following equations yield Assertion~(2):
$$
(Q_{11}):\ 0=-1-C_{11}{}^1,\quad(Q_{12}):\ 0=-C_{12}{}^1,\quad(Q_{22}):\ 0=-C_{22}{}^1\,.
$$

\smallbreak\noindent{\bf Case 3.} Suppose $E(0,\nabla)=\operatorname{Span}\{(x^1)^\alpha,(x^1)^\beta\}$ for $\alpha\ne\beta$.
Since $1\in E(0,\nabla)$, we can take $\beta=0$ and $\alpha\ne0$.
The following equations yield Assertion~(3):
\medbreak\noindent\hfill
$(Q_{11}):\ 0=\alpha(-1+\alpha-C_{11}{}^1),\  (Q_{12}):\ 0=-\alpha C_{12}{}^1,\  (Q_{22}):\ 0=-\alpha C_{22}{}^1$.
\hfill\vphantom{.}\qed

\subsection{The proof of Theorem~\ref{T1.15}: $\mu=-1$}\label{S6.3}
Let  $\mathcal{M}$ be a Type $\mathcal{B}$ surface. By Theorem~\ref{T1.5}, $\dim\{E(-1,\nabla)\}\neq 2,3$.
Suppose  $\dim\{E(-1,\nabla)\}=1$. By Theorem~\ref{T1.12},
$f(x^1,x^2)=(x^1)^\alpha$ for some $\alpha$. As in the proof of Theorem~\ref{T1.13}, we distinguish cases. We clear denominators.
\smallbreak\noindent{\bf Case 1.} Suppose $C_{22}{}^1=0$.
\begin{eqnarray*}
&&(Q_{11}):\ 0=\alpha ^2-\alpha  (C_{11}{}^1+1)+C_{12}{}^2 (C_{11}{}^1-C_{12}{}^2+1)+C_{11}{}^2 (C_{22}{}^2-C_{12}{}^1),\\
&&(Q_{12}):\ 0=C_{12}{}^1 (-2 \alpha +2 C_{12}{}^2-1)+C_{22}{}^2,\ \ 
(Q_{22}):\ 0=C_{12}{}^1 (C_{22}{}^2-C_{12}{}^1)\,.
\end{eqnarray*}
If $C_{12}{}^1=0$, then Equation
$(Q_{12})$ implies $C_{22}{}^2=0$. Thus $C_{12}{}^1=C_{22}{}^1=C_{22}{}^2=0$ so by Theorem~\ref{T1.13}, $\mathcal{M}$
is projectively flat. This is false. Consequently, by Equation $(Q_{22})$, we have 
$C_{12}{}^1=C_{22}{}^2\ne0$. We set $\alpha=C_{12}{}^2\ne0$ to satisfy equations and thereby obtain the structure of Assertion~(1).
\smallbreak\noindent{\bf Case 2.} Suppose $C_{22}{}^1\ne0$. We can rescale so $C_{22}{}^1=\varepsilon=\pm1$ and change coordinates to ensure $C_{12}{}^1=0$. 
We obtain
\medbreak\quad$(Q_{11}):\ 0=\alpha ^2-(C_{11}{}^1+1) \alpha -(C_{12}{}^2)^2+C_{11}{}^1 C_{12}{}^2+C_{12}{}^2+C_{11}{}^2 C_{22}{}^2${ ,}
\smallbreak\quad$(Q_{12}):\ 0=C_{22}{}^2-2 C_{11}{}^2 \epsilon,\qquad
(Q_{22}):\ 0=\alpha -C_{11}{}^1+C_{12}{}^2+1$.
\medbreak\noindent We set $C_{22}{}^2=2C_{11}{}^2\varepsilon$ and obtain
$0=1+\alpha-C_{11}{}^1+C_{12}{}^2$. We set $\alpha=C_{11}{}^1-C_{12}{}^2-1$. The final Equation then yields
$C_{11}{}^1=1+2C_{12}{}^2+\varepsilon(C_{11}{}^2)^2$. This is the structure of Assertion~(2).\qed

\subsection{The proof of Theorem~\ref{T1.17}}\label{S6.4}
Let $\mu\notin\{0,-1\}$. Suppose $\mathcal{M}$ is a Type~$\mathcal{B}$ surface with $\rho_{s,\nabla}\ne0$ and which is not also
of Type~$\mathcal{A}$.
Suppose $\dim\{E(\mu,\nabla)\}\ge1$; $\dim\{E(\mu,\nabla)\}\le2$ by Theorem~\ref{T1.5}. By Theorem~\ref{T1.12},
$(x^1)^\alpha\in E(\mu,\nabla)$. As in the proof of Theorem~\ref{T1.15}, we distinguish cases.

Suppose first that $C_{22}{}^1=0$. Equation~$(Q_{22})$ shows $0=C_{12}{}^1(C_{12}{}^1-C_{22}{}^2)\mu$.
Since $\mu\ne0$, either $C_{12}{}^1$ or $C_{12}{}^1=C_{22}{}^2$. If $C_{12}{}^1=0$, we use Equation~$(Q_{12})$ to see
$0=-C_{22}{}^2\mu$ so $C_{22}{}^2=0$. This gives a Type~$\mathcal{A}$ structure contrary to our assumption. 
We therefore obtain $C_{12}{}^1=C_{22}{}^2$. Equation~$(Q_{12})$ shows $0=-C_{12}{}^1(\alpha+C_{12}{}^2\mu)$.
Thus, $\alpha=-C_{12}{}^2\mu$. Equation~$(Q_{11})$ shows $ 0=(C_{12}{}^2)^2\mu(1+\mu)$.
Since $\mu\notin\{0,-1\}$, $C_{12}{}^2=0$. This implies $\alpha=0$ so $1\in E(\mu,\nabla)$ and $\rho_s=0$ contrary to our assumption.

We therefore have $C_{22}{}^1\ne0$. We can rescale the coordinates so $C_{22}{}^1=\varepsilon=\pm1$. We can then
make a shear so $C_{12}{}^1=0$. We substitute these relations to obtain:
\medbreak\qquad $(Q_{11}):\ 0=\alpha ^2-\alpha  (C_{11}{}^1+1)-\mu  \left(C_{11}{}^1 C_{12}{}^2+C_{11}{}^2 C_{22}{}^2-(C_{12}{}^2)^2+C_{12}{}^2\right)$,
\smallbreak\qquad $(Q_{12}):\ 0=C_{11}{}^2 \mu  \epsilon -\frac{C_{22}{}^2 \mu }{2}$,\qquad
$(Q_{22}):\ 0=\epsilon  (\mu  (-C_{11}{}^1+C_{12}{}^2+1)-\alpha )$.
\medbreak\noindent 
Since $\mu\ne0$, we have $C_{22}{}^2=2\varepsilon C_{11}{}^2$. We have:
\begin{equation}\label{E6.a}
C_{22}{}^1=\varepsilon,\quad C_{12}{}^1=0,\quad \alpha=\mu(1+C_{12}{}^2-C_{11}{}^1)\,.
\end{equation}
The only remaining Equation is
\medbreak\quad$(Q_{11}):\ 0=\mu\{(C_{11}{}^1)^2 (\mu +1)-2 C_{11}{}^1 (C_{12}{}^2 \mu +C_{12}{}^2+\mu )-2 (C_{11}{}^2)^2 \epsilon$
\smallbreak\qquad\qquad\qquad$+(C_{12}{}^2)^2 (\mu +1)+2 C_{12}{}^2 (\mu -1)+\mu -1\}$.
\medbreak\noindent
Since $\mu\ne0$, we can solve Equation~($Q_{11}$) for $\mu$ to complete the proof of the first portion of Theorem~\ref{T1.17}.

We now suppose $\dim\{E(\mu,\nabla)\}=2$. We examine the possibilities
of Theorem~\ref{T1.13} seriatim.
\subsection*{Case 1} $E(\mu,\nabla)=\operatorname{Span}\{(x^1)^\alpha,(x^1)^\alpha(cx^1+x^2)\}$. Let $f=(x^1)^\alpha(cx^1+x^2)$. We have
$\alpha=(1-C_{11}{}^1+C_{12}{}^2)\mu$. Equation~($Q_{22}$) shows $0=-(c+2C_{11}{}^2)$
so  $c=-2 C_{11}{}^2$. After clearing denominators, we have
\medbreak\quad$(Q_{11}):\ 0=C_{11}{}^2 \{-2 (C_{11}{}^1)^2+C_{11}{}^1 (6 C_{12}{}^2-3)+8 \varepsilon (C_{11}{}^2)^2$
\smallbreak\qquad\qquad\qquad$-4 (C_{12}{}^2)^2+9 C_{12}{}^2+5\}$,
\smallbreak\quad$(Q_{12}):\ 0=(C_{11}{}^1)^2-3 C_{11}{}^1 C_{12}{}^2-2\varepsilon (C_{11}{}^2)^2+2 (C_{12}{}^2)^2-C_{12}{}^2-1$.

\subsection*{Case 1a}
If $C_{11}{}^2=0$, Equation~($Q_{11}$) is trivial and we obtain
\medbreak\quad$(Q_{12}):\ 0=(C_{11}{}^1)^2-3 C_{11}{}^1 C_{12}{}^2+2 (C_{12}{}^2)^2-C_{12}{}^2-1$
\smallbreak\quad\qquad\quad\phantom{aaa}$=(C_{11}{}^1-2 C_{12}{}^2-1) (C_{11}{}^1-C_{12}{}^2+1)$.
\medbreak\noindent
If $C_{11}{}^1=1+2C_{12}{}^2$, then $\mu=-1$ which is false. If $C_{11}{}^1=C_{12}{}^2-1$, then we obtain the structure in Assertion~(1).
\subsection*{Case 1b} Suppose $C_{11}{}^2\ne0$, we may divide the first equation~by $C_{11}{}^2$ to see
\medbreak\quad$(\tilde Q_{11}):\ 0=-2 (C_{11}{}^1)^2+C_{11}{}^1 (6 C_{12}{}^2-3)+8\varepsilon (C_{11}{}^2)^2-4 (C_{12}{}^2)^2+9 C_{12}{}^2+5$.
\medbreak\noindent We compute that:
\medbreak\quad$(\tilde Q_{11})+4(Q_{12}):\ 0=2 (C_{11}{}^1)^2-3 C_{11}{}^1 (2 C_{12}{}^2+1)+4 (C_{12}{}^2)^2+5 C_{12}{}^2+1$
\smallbreak\quad\qquad\qquad\qquad\qquad\phantom{..}$=(2 C_{11}{}^1-4 C_{12}{}^2-1) (C_{11}{}^1-C_{12}{}^2-1)$.
\medbreak\noindent Since $(C_{11}{}^1-C_{12}{}^2-1)\ne0$, we obtain
$2 C_{11}{}^1-4 C_{12}{}^2-1=0$. There is then a single remaining relation:
$0=8 (C_{11}{}^2)^2 \epsilon +2 C_{12}{}^2+3$. We solve this for $C_{12}{}^2$ to obtain the structure of
Assertion~(2).
\subsection*{Case 2} $E(\mu,\nabla)=\operatorname{Span}\{(x^1)^\alpha,(x^1)^\alpha\log(x^1)\}$. Evaluating Equation~($Q_{22}$) at $x^1=1$ yields
$\varepsilon=0$ which is impossible. Therefore this case does not arise.
\subsection*{Case 3} $E(\mu,\nabla)=\operatorname{Span}\{(x^1)^\alpha,(x^1)^\beta\}$ for $\alpha\ne\beta$. 
Since $\alpha$ is determined by Equation~(\ref{E6.a}), this case does not arise.

\end{document}